\documentclass[12pt]{article}
\usepackage{amssymb}
\usepackage{amsmath}
\usepackage{latexsym}
\usepackage{showkeys}

\newtheorem{theorem}{Theorem}[section]

\newtheorem{proposition}[theorem]{Proposition}
\newtheorem{corollary}[theorem]{Corollary}

\def\eps{{\varepsilon}}

\def\msn{\medskip\noindent}
\def\nrv{net reproductive rate}

\def\rarr{\rightarrow}
\newcommand{\bequ}{\begin{equation}}
\newcommand{\eequ}{\end{equation}}
\newcommand{\bproof}{\underline{\em Proof}:\ }
\newcommand{\eproof}{{\hfill $\Box$}\\ }

\begin{document}

\title{Applications of Perron-Frobenius Theory\\ to Population
Dynamics}

\author{Chi-Kwong Li\thanks{Research partially supported by an NSF grant}
\\Department of Mathematics\\ College of William and Mary\\ P.O.
Box 8795\\
 Williamsburg, VA 23187-8795, USA,\\E-mail: ckli@math.wm.edu
 \and Hans Schneider\thanks{Research supported by the U.S. Social Security
Administration and the Wisconsin Department of Employee Trust
Funds}\\ Department of Mathematics\\ University of Wisconsin\\ 480
Lincoln Drive\\ Madison, WI 53706, USA\\E-mail:
hans@math.wisc.edu}

\date{23 April 2001}

\maketitle

\abstract By the use of Perron-Frobenius theory, simple proofs are
given of the Fundamental Theorem of Demography and of a theorem of
Cushing and Yicang on the net reproductive rate occurring in
matrix models of population dynamics. The latter result is further
refined with some additional nonnegative matrix theory. When the
fertility matrix is scaled by the net reproductive rate, the
growth rate of the model is $1$. More generally, we show how to
achieve a given growth rate for the model by scaling the fertility
matrix. Demographic interpretations of the results are given.
\section{Introduction}
A \em standard matrix model of population dynamics \rm is given by
a sequence of nonnegative vectors $x_0, x_1, \ldots $ of fixed
length $n$ defined by
\begin{equation}\label{equ1}
x_k = P x_{k-1}, \quad k = 1, 2, \dots,
\end{equation}
for a given nonzero $x_0$ where $P$ is an $n \times n$ matrix with
nonnegative entries. As usual, we assume that
\begin{equation} \label{equ2}  P = T+ F
\end{equation}
where $T$ and $F$ are nonnegative nonzero matrices such that all
the column sums of $T$ are not larger than one:
\begin{equation} \label{equ3}
\Sigma_{i=1,\ldots,n}\ t_{ij}\leq 1, \quad j = 1,\ldots,n.
\end{equation}
The $m$th entry of the vector $x_k$ represents the population in
the $m$th class at the time epoch $k$. The matrix $P = T + F$ is
known as the {\it projection matrix}; the matrices $T$ and $F$ are
known as the {\it transition matrix} and the {\it fertility
matrix} resp., so that the $(i,j)$ entry of $T$ represents the
fraction of the individuals in the $j$th class that will survive
and move to the $i$th class in a unit interval of time, and the
$(i,j)$ entry of $F$ represents the number of newborns in the
$i$th class that descend from one individual in the $j$th class in
a unit interval. Because of these demographic interpretations, one
sees why the column sums of $T$ are always less than or equal to
1. The special case of this matrix model for Leslie matrices,
introduced by Bernardelli \cite{Ber}, Lewis \cite{Lew} and {Leslie
\cite {L1}, has a long history
  \footnote{ \cite{L1} acknowledges the contribution of
\cite{Ber} and \cite{L2} mentions \cite{Lew}, but \cite{Ber} and
\cite{Lew} do not refer to each other, which is extraordinary as
both authors were at University College, Rangoon, Burma, in {1940
- 41} when their papers were submitted.}}. We concentrate on
properties of the general case introduced above.

\msn The following argument shows that a further assumption is
warranted for our matrix model. The {\it spectral radius}
$\rho(T)$ of a matrix $T$ is the maximum of the moduli of its
eigenvalues. Since we wish to exclude the possibility of an
immortal population, we shall always assume in the rest of the
paper that $\lim_{k \rarr \infty}T^kx_0=0$ for all initial
populations (nonnegative vectors) $x_0$. By examining the effect
of $T^k$ on the standard basis of unit vectors, this condition may
be shown to be equivalent to $\lim_{k \rarr \infty}T^k=0$. In
turn, it is known that this condition is equivalent to
\begin{equation} \label{mort}
\rho(T) < 1 ,
\end{equation}
 \cite[Theorem 3.5]{W2}, \cite[Theorem 5.612]{HJ}.
We shall thus always assume that our model satisfies (\ref{mort})
together with (\ref{equ1})and (\ref{equ2}); the assumption
(\ref{equ3}) is not used in our proofs. We remark that
(\ref{equ3}) implies the weaker condition $\rho(T) \leq 1$, but
that (\ref{mort}) does not imply (\ref{equ3}) as is easily shown
by examples.

\msn Under the assumption (\ref{mort}) we have
\begin{equation} \label{defq}
(I-T)^{-1} \ = \ I+T + T^2 + \ldots,
\end{equation}
e.g.\ \cite[Theorem 6.1]{W2}, \cite[Corollary 5.6.16]{HJ}. Let $Q
= F(I-T)^{-1}$. Then we have
 $$ Qx_0 = Fx_0 + FTx_0 + FT^2x_0 + \ldots,$$
which represents the distribution with respect to state-at-birth
of all newborn descendants accumulated during the entire lifespan
of the population $x_0$, see \cite[p.8]{C} and \cite[p.71]{DH},
where the matrix $Q$ is called the ``next generation matrix''.
Following \cite[p. 126]{Cas2}, we call the spectral radius
$\rho(Q)$ of $Q$ the {\em net reproductive rate} of the model, see
\cite{CY} and \cite[p.8]{C}, for a definition in this generality.
As usual it is denoted by $R_0$. Other names for $R_0$ are ``net
reproductive number'' (\cite{C}) and ``basic reproduction ratio"
(\cite{DH}).

\msn The purpose of this note is to explore some applications of
the Perron-Frobenius theory to results in Population Dynamics;
particularly we wish to study the role of the net reproductive
rate. In Section 2 we state needed standard results from this
theory which will be used throughout the paper. We first apply the
theory to obtain the Fundamental Theorem of Demography for a model
with a primitive projection matrix. This approach is well-known,
see for example \cite{Pol}, \cite{I}, or \cite{Cas2}. In Section
3, we consider the net reproductive rate. We state and give a
short proof of a somewhat stronger form of an interesting
comparison theorem on the net reproductive rate due to Cushing and
Yicang, see \cite[Theorem 3]{CY} and \cite[Theorem 1.1.3]{C}. The
theorem is further refined in Section 4 using a result from graph
theoretic Perron-Frobenius theory. For a model with an irreducible
projection matrix, the net reproductive rate may be viewed as a
factor producing a model with growth rate $1$ when one scales only
the fertility matrix by this constant factor. More generally,
given a positive $s$ subject to one restriction, we determine a
$q(s)$ as a function of $s$ such that the growth rate of the model
is $s$ when the fertility matrix is scaled by $q(s)$. In our last
section, we give demographic interpretations of our results.

\section{Perron-Frobenius Theory and the Fundamental Theorem of Population
Demography}  \label{PF0}

In this section we give the sketch of a proof of the fundamental
theorem of demography, see \cite[Theorem 1.1.2]{C} or
\cite[p.86]{Cas2} using the Perron-Frobenius theory of nonnegative
matrices. This theory is also needed for applications in
subsequent sections.

\msn A nonnegative matrix is {\em irreducible} if it is not the $1
\times 1$ zero matrix and it is not permutationally similar to a
matrix of the form
 $$\begin{pmatrix}
 A_{11} & A_{12} \\ 0 & A_{22}
 \end{pmatrix} $$
for nontrivial square matrices $A_{11}$ and $A_{22}$.  The
celebrated Perron-Frobenius theorem for irreducible matrices  can
be stated as follows.

\begin{theorem} \label{PF1} Let $P$ be an irreducible nonnegative matrix.
Then

\msn {\rm (a)} The spectral radius $\rho(P)$ of $P$ is positive
and it is an algebraically simple eigenvalue of $P$ with
corresponding left and right positive eigenvector, which are
unique up to scalar multiples.

\msn {\rm (b)} The spectral radius of $P$ is the unique eigenvalue
with a nonnegative eigenvector.

\msn {\rm (c)} The spectral radius of the matrix $P$ increases
(strictly), resp. decreases, if any entry of it increases, resp.
decreases. \eproof \end{theorem}
 One may see \cite[Theorem
10.7]{W2}, \cite[Theorem 8.4.4, Problem 15, p.515]{HJ},
\cite[Theorem 2.1]{V}, \cite[Theorem 2, p.53]{G} for proofs.

\msn Following a common practice in matrix literature, we call the
spectral radius $\rho(P)$ of a nonnegative matrix $P$ the {\em
Perron root} of $P$. In some fields of mathematics this term is
defined as the reciprocal of the radius of convergence of the
power series $\Sigma_0^\infty  z^r P^r$ but by a classical result
this is precisely $\rho(P)$ when $P$ is a finite complex matrix,
see \cite[Theorem 5.5]{W2}, or see the more general Theorem by
Hensel \cite[Theorem 5.4]{W2}. We also call a nonnegative left or
right eigenvector corresponding to the Perron root of a
nonnegative matrix a {\em Perron vector}. We usually denote the
Perron root of a projection matrix  $P$ by $r$ and left and right
Perron vectors by $v^t$ and $u$ respectively.

\msn Since the eigenvalues and eigenvectors (with proper
normalization) are continuous in the elements of a matrix and
every nonnegative matrix is a limit of a sequence of irreducible
nonnegative matrices, one immediately deduces:
\begin{corollary} \label{PF2} Let $P$ be a nonnegative matrix.
Then

\msn {\rm (a)} The spectral radius $\rho(P)$ of $P$ is an
eigenvalue of $P$ with a corresponding nonnegative Perron vector.

\msn {\rm (b)} The spectral radius of the matrix $P$ does not
decrease, resp. increase, if any entry of it increases, resp.
decreases. \eproof
\end{corollary}

\msn An irreducible nonnegative matrix $P$ is {\it primitive} if
there is only one eigenvalue of $P$ that attains the modulus
$\rho(P)$. In this case, the Perron root is properly called the
{\em dominant eigenvalue} of $A$, which is the usual term in the
demographic literature (even under less restrictive conditions),
e.g. \cite[p. 83]{Cas2}. A nonnegative matrix $P$ is primitive if
and only if $P^k$ is positive for some positive integer $k$,
\cite[Theorem 8.5.2]{HJ}. This shows that primitivity (like
irreducibility) depends only on the {\em pattern} of a nonnegative
matrix, i.e. if $P$ is primitive (irreducible) then every matrix
that has positive entries in exactly the same positions is also
primitive (irreducible). Applying a standard result on the
convergence of powers of matrices \cite[Theorem 3.5]{W2}, one may
derive from Theorem \ref{PF1} the following well known result, for
closely related results see for example \cite[p.86]{Cas2}, or in
the special case of Leslie matrices, \cite[Lemma 4.5.1]{Pol} or
\cite{I}, where proofs are given in the same spirit as ours.
\begin{theorem} \label{PF3} Let $P$ be a primitive nonnegative matrix
with spectral radius $\rho(P)=r$ and left and right Perron vectors
$v^t$ and $u$ such that $v^tu = 1$. Then
 \begin{equation} \label{lim}
  \lim_{k \rightarrow \infty}(P/r)^k =  uv^t.
 \end{equation}
\eproof
\end{theorem}
The following application of Theorem \ref{PF3} is called the
{\em fundamental theorem of demography} in \cite[Theorem
1.1.2]{C}, where a different proof is given.
\begin{theorem} \label{PF4} Let $P$ be the projection matrix of a standard
population model $x_k, \ k = 0,1, \ldots$, given by {\rm
(\ref{equ1})}. Suppose that $P$ is primitive with spectral radius
$\rho(P)=r$ and has left and right Perron vectors $v^t$ and $u$
resp. normalized so that $v^tu = 1$. Then
 \begin{equation*}
\lim_{k \rightarrow \infty}x_k/r^k =  (v^tx_0)u.
 \end{equation*}
Consequently, if $|w|$ denotes the sum of entries of the vector $w$, so
that $|x_k|$ will denote the total population at time $k$ in
the population model, then
$$\lim_{k \rightarrow \infty} |x_k|
= \left\{ \begin{matrix}
 0 \hfill & {\rm if }\ r < 1,\\
 |(v^tx_0)u| \hfill  & {\rm if }\ r = 1, \\
 \infty \hfill  &  {\rm if }\  r > 1.
\end{matrix} \right.
$$  \eproof
\end{theorem}

\msn  Mathematically, a {\em population} is a vector that is
nonnegative and nonzero. As is usual in the demographic
literature, we call a population $x$ {\em stable} (for a given
matrix model), though this terminology is not consistent with
definitions of stability in other parts of mathematics. If, for
some positive $r,\ Px = rx$, and we call $x$ a {\em stationary
population} if $Px = x$ . We call a population {\em eventually
stable} if $\lim_{k \rightarrow \infty}x_k/r^k$ exists and is
nonzero, and we call the population {\em eventually stationary} if
$\lim_{k \rightarrow \infty}x_k$ exists and is nonzero. The
spectral radius $\rho(P) = r$ is called the {\em growth rate} of
the model, which for primitive $P$ is justified by Theorem
\ref{PF4}.

\msn Suppose that $P$ is primitive. It follows immediately from
the fundamental theorem that, whatever the initial population
$x_0$, the number of individuals in $x_k$ grows to infinity if $r
> 1$, shrinks to $0$ if $r < 1$, and remains finite if $r=1$.
Furthermore, in  all cases, there exists a unique stable
population (except for a constant factor), which is a stationary
population if $r=1$. Furthermore, in this case, the fundamental
theorem shows that every population is eventually stationary.

\msn The assumption that $P$ is primitive cannot be omitted from
this last remark, for the conclusion (\ref{lim}) of Theorem
\ref{PF4} depends strongly on this assumption. If one merely
assumes that $P$ is  irreducible  one may show that $x_k/r^k, k =
1,2 \ldots $ is bounded above.  For reducible nonnegative $P$ the
description of the possible limiting behavior of $P^kx_0$ is quite
complicated, see \cite{FS} for applications of  graph theoretic
concepts to this problem. An example is given near the end of this
article where $r = 1$, but all populations except for the stable
populations grow to be infinitely large.

\msn For primitive $P$, the reciprocal of the growth factor $r$ of
$P$ may also be interpreted as a factor for stationarity for the
model, viz. if both $T$ and $F$ are scaled by the same factor
$1/r$ then the resultant model with matrix $P' = (T+F)/r$ has the
property that every population is eventually stationary. We are
however particularly interested in scaling the fertility matrix
$F$ without scaling the transition matrix $T$ so that in the
resultant model every population is eventually stationary, and
this leads naturally to the considerations in the rest of this
paper.

\section{The net reproductive rate}
The main result of this section is a somewhat stronger form of a
comparison theorem due to Cushing and Yicang, \cite[Theorem3]{CY},
see also \cite[Theorem 1.1.3]{C}. Using standard results of
Perron-Frobenius theory reviewed in Section \ref{PF0}, we give a
very short and simple proof of this theorem stated below as
Theorem \ref{comp0}. In the next section we show that this theorem
can be further refined and generalized using some more nonnegative
matrix theory.
 \begin{theorem} \label{comp0}  Suppose a standard matrix model
of population dynamics satisfies  {\rm (\ref{equ1})} and {\rm
(\ref{equ2})}, and assume that the projection matrix $P = T+F$ is
irreducible where  $T$ is nonzero and satisfies {\rm
(\ref{mort})}. Denote the growth factor $\rho(P)$ by $r$ and the
net reproductive rate $\rho(Q)$, where $Q = F(I-T)^{-1}$, by
$R_0$. Suppose that $R_0
> 0$. Then
\begin{equation} \label{stable}
\rho(T+F/R_0) = 1,
\end{equation}
and  one of the following holds:
\begin{equation}  \label{comp1}
 r=R_0=1, \qquad {\rm or} \qquad  1 < r < R_0, \qquad {\rm or} \qquad 0<R_0<r<1.
\end{equation}
\end{theorem}
\bproof Since $\rho(T) < 1$, by (\ref{defq}) the matrix
$(I-T)^{-1}$ is nonnegative and hence so is $Q = F(I-T)^{-1}$.
Clearly $F \neq 0$ since $\rho(Q) > 0$.

\msn To prove the equality (\ref{stable}), note that by Corollary
\ref{PF2}, there exists a nonnegative left eigenvector $y^t$ of
$F(I-T)^{-1}$ corresponding to the eigenvalue $R_0$, i.e.
$y^tF(I-T)^{-1} = R_0 y^t $. Then $y^tF = R_0y^t(I-T)$ and hence
$y^t(T+F/R_0) = y^t$. Since $R_0>0$, the matrix $T+F/R_0$ is
irreducible and hence it follows by Theorem \ref{PF1}(b) that
$\rho(T +F/R_0) = 1$.

\msn To prove that one of the conditions in (\ref{comp1}) holds we
consider three cases.

\medskip\noindent
(i) \ If $R_0 = 1$, then $1 = \rho(T+F) = r$.

\medskip\noindent
(ii) \ If $R_0 > 1$, then $$T+F/R_0 \leq T+F  \leq R_0 T+F$$
 with equalities excluded since $F$ is nonzero. Hence  by Theorem
 \ref{PF1}(c)
  $$ 1 = \rho(T+F/R_0) < \rho(T + F) = r  < \rho(R_0 T+F) =R_0.$$

\msn (iii) \ If $0 < R_0 < 1$, then again by Theorem \ref{PF1}(c),
 $$ 1 = \rho(T+F/R_0) > \rho(T+F) = r  >  \rho(R_0 T+F) = R_0.$$
\eproof
 The special case of Theorem \ref{comp0} when $T$ is
strictly lower triangular and F is upper triangular is known in
numerical linear algebra as the Stein-Rosenberg Theorem,see
\cite[pp.68-70]{V} for a proof to which our proof of the more
general result is somewhat similar. See \cite{RV} for a result
close to Theorem \ref{comp0}.

\begin{corollary} \label{forget}
Under the hypotheses of Theorem {\rm \ref{comp0}}, consider the
modified model given by $\tilde P = T + F/R_0$ with left and right
Perron vectors $\tilde v^t$ and $\tilde u$ such that $\tilde v^t
\tilde u = 1$. If $P$ is primitive, then for every initial
population $x_0$ we have
 \begin{equation*}
 \lim_{k \rarr \infty}(T+F/R_0)^kx_0 = (\tilde v^t x_0) \tilde u .
 \end{equation*}
\end{corollary}
\bproof This follows immediately by Theorem \ref{comp0}(a) and
Theorem \ref{PF3}, since $\tilde P$ is also primitive. \eproof

\msn There is a corresponding theorem for general nonnegative
matrices.
 \begin{theorem} \label{compg}  Suppose a standard matrix model
of population dynamics satisfies  $(\ref{equ1})$ and
$(\ref{equ2})$, and the transition matrix $T$ satisfies
$(\ref{mort})$. Denote the growth factor $\rho(P)$ by $r$ and the
net reproductive rate $\rho(Q)$, where $Q = F(I-T)^{-1}$, by
$R_0$. Then one of the following holds:
\begin{equation}  \label{compg1}
 r=R_0=1,\qquad {\rm or} \qquad 1 < r \leq R_0, \qquad {\rm or} \qquad  0 \leq R_0
  \leq r  <1.
 \end{equation}
 If $R_0 > 0$, then
 \begin{equation} \label{compg0}
 \rho(T+F/R_0) = 1.
\end{equation} \eproof
\end{theorem}
Most of the derivation of this theorem by continuity from Theorem
\ref{comp0} is straightforward and therefore omitted. However, we
shall show why $r=1$ implies that $R_0 = 1$. For positive $\eps$,
let $F(\eps) = F + \eps E$, where $E$ is a matrix of the
appropriate size all of whose entries are $1$. Let $P(\eps) = T +
F(\eps)$. Since $P(\eps)$ is irreducible, we have $\rho(P(\eps)) >
\rho(P) = 1$, and hence by Theorem \ref{comp0},
$\rho(F(\eps)(I-T)^{-1}) > 1$. If $r = 1$, letting $\eps$ tend to
$0$ we obtain $R_0 \geq 1$. Now consider $P'(\eps) =
P(\eps)/\rho(P(2\eps))$. Since $\rho(P(2\eps)) > \rho(P(\eps))$ we
have $\rho(P'(\eps)) < 1$ and hence
$\rho(F/\rho(P(2\eps))(I-T/\rho(P(2\eps)^{-1}))) < 1$, again by
Theorem \ref{comp0}. Letting $\eps$ tend to $0$ we now obtain $R_0
\leq r = 1$. It follows that $R_0 = 1$. We observe that Theorem
\ref{compg} may also be deduced by means of \cite[Theorem
4.5]{Sch}.

\section{The refined stability and comparison theorem}
\label{refined}
 In this section we show that the hypothesis
$R_0 = \rho(Q)>0$ in Theorem \ref{comp0} actually follows from the
remaining assumptions thus allowing us to state a refined version
of this theorem. Our proof depends on the following proposition,
which is a restatement of \cite[Lemma 3.4]{Sch}, see also
\cite{Sz} and \cite{Z}.
\begin{proposition} \label{split1} Let $T$  and $F$ be nonnegative
matrices with $\rho(T) < 1$ and $F \neq 0$. Suppose
$T+F$ is irreducible and $Q = F(I-T)^{-1}$. Then, after a permutation
similarity,
\begin{equation} \label{split0}
Q = \begin{pmatrix}
        Q_{11}  &  Q_{12} \\
        0       &  0
\end{pmatrix},
\end{equation}
where $Q_{11}$ is a nontrivial irreducible nonnegative matrix,
$Q_{12}$ is a nonnegative matrix every column of which has a
positive entry, and the $0$ rows of $Q$ correspond to the $0$ rows
of $F$, if any.  \eproof
\end{proposition}
Our stability and comparison theorem may now be stated as:
\begin{theorem} \label{compfin}
Suppose a standard matrix model of population dynamics satisfies
{\rm (\ref{equ1})} and {\rm (\ref{equ2})}, and assume that the
projection matrix $P = T+F$ is irreducible with $T$ satisfying
{\rm (\ref{mort})} and $F \neq 0$. Denote the growth rate
$\rho(P)$ by $r$ and the net reproductive rate $\rho(Q)$, where $Q
= F(I-T)^{-1}$, by $R_0$. Then equation $(\ref{stable})$ holds and
so does one of the conditions in $(\ref{comp1})$. Furthermore, the
matrix $Q$ is irreducible if and only if every row of $F$ contains
a positive element.
\end{theorem}
\bproof We have $\rho(Q) = \rho(Q_{11})$ and since $Q_{11}$ is
irreducible, $\rho(Q_{11})>0$ by Theorem \ref{PF1}. The first
conclusion follows from Theorem \ref{comp0} and the second from
Proposition \ref{split1} (which was also observed in
\cite{Z}).\eproof

\msn Seneta \cite[p.42]{Sen} gives a proof of the Stein-Rosenberg
theorem which may be adapted to show $R_0 > 0 $ under the
hypotheses of Theorem \ref{compfin}. The key ingredient is the
observation that $(I-T)^{-1}$ and $(I-T/r)^{-1}$ have the same
zero--nonzero pattern.

\msn In contrast to the situation for irreducible $P$, one may
construct examples of a  reducible projection matrix $P$ with a
nonzero fertility matrix $F$ such that the corresponding net
reproductive rate  $ R_0=0$. In this connection,  we shall prove
the following theorem.
\begin{theorem} \label{zero}
Let $P, T$ and $F$ satisfy the hypotheses of Theorem {\rm
\ref{compg}}. Then the net reproductive rate $R_0 > 0$ if and only
if for some $a>0,\ \rho(T+aF) > \rho(T)$.
 \end{theorem}
\bproof First suppose that $R_0 > 0$. Then, for sufficiently large
positive $a$, we have $aR_0 > 1$. But this is the net reproductive
rate of the projection matrix $T+aF$ and hence, by Theorem
\ref{compg}, $\rho(T+aF) > 1 > \rho(T)$.

\msn Conversely, suppose $\rho(T+aF) > \rho(T)$ for some positive
number $a$.  Let $\det(\lambda I - (T+aF)) = \sum_{k=1}^n
\lambda^k f_k(a)$. Using the usual determinantal expansions, we
see that the $f_k(a),\ k =1, \ldots, n$,  are polynomials in  $a$.
But they are also signed sums  of the $k$-th elementary symmetric
functions of the eigenvalues of the matrix $T+aF$.  Since
$\rho(T+aF) > \rho(T)$, there exists $k$ such that $f_k(a)$ is not
a constant polynomial, and hence, for this $k$,  $|f_k(a)|$ is
unbounded as $a$ goes to infinity. Hence at least one eigenvalue
cannot be bounded in $a$, and it follows that $\rho(T+aF)$ is
unbounded, and thus there exists $a$ such that $\rho(T+aF) > 1$.
Again applying Theorem \ref{compg} to $T+aF$, we obtain $aR_0 > 1$
and hence $R_0 > 0$. \eproof

\msn Our proof of Theorem \ref{zero} also shows that if $R_0$ is
positive, then $F$ may be scaled to achieve an arbitrarily large
growth rate in the general nonnegative case. Returning to an
irreducible projection matrix $P = T+F$ with $F \neq 0$ we observe
that Theorem \ref{zero} provides a second proof that $R_0 > 0$,
since in this case $\rho(T + aF) > \rho(T)$ for all $a > 0$.

\msn We now generalize Theorem \ref {compfin}. If $F+T$ is
irreducible with Perron root $r$ and positive left Perron vector
$z^t$ then it easily follows that $z^tF(I-T/r)^{-1} = rz^t$ and
since $z^t > 0$ we deduce that $\rho(F(I-T/r)^{-1}) = r$, see e.g.
\cite[Cor. 8.1.30]{HJ}. Thus Theorem \ref{compfin} is the special
case $s =r,\ q(s) = 1$ of the result which now follows.
 \begin{theorem} \label{genthm}
 Let $P,T$ and $F$ satisfy the conditions of Theorem {\rm \ref{compfin}}.
 For $s > \rho(T)$ define
  \begin{equation} \label{assign}
  q(s) = \rho(F(I-T/s)^{-1})/s .
  \end{equation}
Then $q(s) > 0$. Let $P(s) = T +F/q(s)$. Then its growth rate,
$\rho(P(s))$, is $s$, and its net reproductive rate is
 $$R_0(s)= R_0/q(s).$$
Further, one of the following holds:
\begin{equation}  \label{comp1a}
1=s= R_0(s), \qquad {\rm or} \qquad  1 < s < R_0(s), \qquad {\rm
or}
 \qquad 0<R_0(s)<s<1.
\end{equation}
\end{theorem}
\bproof We observe that $Q(s) = F(I-T/s)^{-1}$ is nonnegative
since $\rho(T/s) < 1$ and hence $q(s)> 0$ by Theorem
\ref{compfin},s as $F + T/s$ is irreducible. Let $z^t$ be the left
Perron vector of $Q(s)$. Thus $z^tQ(s) = z^tF(I-T/s)^{-1} = s
q(s)z^t$. An easy computation now yields
 $$z^t(T+F/q(s)) = sz^t $$
and hence $s$ is the Perron root of $T+F/q(s)$ as asserted in the
theorem. The corresponding net reproductive rate is $R_0(s) =
\rho((F/q(s))(I-T)^{-1}) = R_0/q(s)$. Then (\ref{comp1a}) follows
by Theorem \ref{compfin}.  \eproof

\msn We note that the inequalities (\ref{comp1a}) are strict and
are equivalent to $\rho(F(I-T)^{-1}) < \rho(F(I-T/s)^{-1})$ if
 $\rho(T) <s < 1$ and $\rho(F(I-T)^{-1}) > \rho(F(I-T/s)^{-1})$ if $s
> 1$. Further, since $\rho(F(I-T/s)^{-1})$  is a decreasing function of $s$,
it follows that $q(s)$ is a strictly decreasing function of $s$
and that $\lim_{s \rightarrow \infty} q(s) = 0$.

\section{Demographic interpretations}
We now discuss the demographic interpretations of our theorems. As
we shall impose various assumptions on the projection matrix $P$,
we begin each paragraph with an assumption which holds throughout
the paragraph.

\msn Assume that $P$ is a primitive nonnegative matrix. As
previously observed, the reciprocal of $r = \rho(P)$ may be viewed
as a factor for the scaling of the fertility and transition
matrices in order to obtain a model with every initial population
eventually stationary population, viz. $\rho(P/r) = 1$. Similarly,
since $R_0 = \rho(Q) > 0$, (Theorem \ref{compfin}) we may leave
the transition matrix fixed and scale the fertility matrix by the
reciprocal of $R_0$ and in the resultant model every initial
population will tend to a stationary population which is unique
except for a multiplicative constant (Corollary \ref{forget}).
Intuitively, scaling the fertility matrix only should be more
radical than scaling both the fertility and the transition matrix
to achieve the same objective, and that is exactly the content of
Theorems \ref{comp0} and \ref{compg}.

\msn Suppose now that $P$ is irreducible. It is possible to adapt
the previous paragraph to this case by considering scalings to
achieve a growth rate of $1$,  but we shall immediately turn to
interpeting the more general Theorem \ref{genthm}. Given any $s
> \rho(T)$, it is possible to obtain a growth rate of $s$ without
changing the transition matrix $T$ by scaling the fertility matrix
$F$ to $F/q(s)$, where $q(s)$ is given by (\ref{assign}). Since
$q(s)$ is a strictly decreasing function of $s$ for $s > \rho(T)$,
to achieve a higher growth rate we require greater fertility,
which again is intuitively clear. We shall discuss these points
further when we turn to the Leslie model.

\msn Suppose that $P$ is irreducible. By Proposition \ref{split1}
every column of the submatrix $(Q_{11}, Q_{12})$ has a positive
entry, and this implies that every class in the initial population
has descendants during its lifetime. By the last assertion of
Theorem \ref{compfin} the matrix $Q$ is irreducible if and only if
there are newborns in every population class in the case of a
population that has members in each class.

\msn Now assume again that $P$ is irreducible. We turn to
characterizations of $R_0$ which may be obtained from Theorem
\ref{compfin} and from the  following characterization of the
Perron root of an irreducible nonnegative matrix $A$ which is a
reformulation of a well-known characterization due to Wielandt
\cite{W1}, see also \cite[p.65]{G}:
  \begin{eqnarray*}
  \rho(A) & = &
 \max \{s:  Ax \geq sx,\ \hbox{ for some}\ x \geq 0,\ x  \neq 0\}\\
 & = & \min \{s:\  Ax \leq sx,\ \hbox{ for some}\ x \geq 0,\ x  \neq
 0\},
  \end{eqnarray*}
 and the equality is attained in either inequality if and only if $x$ is
the Perron vector of $A$. For this formulation and references to
it, see the commentary following \cite{W1} in Wielandt's
Mathematical Works. Let $z$ denote a nonnegative vector with the
same number of entries as there are columns in  $Q_{11}$, and let
$z_i$ be its $i$th entry. Then there exist indices $i$ and $j$
such that $(Q_{11}z)_i \geq R_0 z_i$ and $(Q_{11}z)_j \leq R_0
z_j$. We may put the case $r = R_0 = 1$ into words thus: For {\em
every} initial population, there must be one class of newborns
that over its lifetime  produces at least as many descendants in
the same class and there must be one class of newborns that over
its lifetime  produces at most as many descendants in the same
class. Furthermore, {\em there exists} a population of newborns
(i.e. corresponding to the nonzero rows of $F$) which over the
course of its lifetime reproduces itself exactly if and only if
there exists a stationary population (but note that this
population of newborns is in general not part of a stationary
population).

\msn Still under the assumption that $P$ is irreducible, it is
easily proved using an additional part of the Perron-Frobenius
theorem, \cite[Theorem 2.20]{BP}, \cite[Corollary 8.4.6 and Remark
8.4.9] {HJ}, \cite[Theorem 2.3]{V}, \cite[Theorem 2, p.53]{G},
that for some positive integer $d$, $P^d$ is the direct sum of
primitive matrices, and hence that there is a population that is
stable in a periodic sense for this model, that is for some
positive integer $d$ and every initial population $x_0$ there
exist populations $w_0, \ldots, w_{d-1}$ (depending on $x_0$) such
that $(\lim_{k \rarr \infty} x_{kd+i})/r^k = w_i, \ i = 0,\ldots,
d-1$. To obtain the previous primitive case, we put $d = 1$.
Another equivalent condition for $\rho(P) = 1$ in the irreducible
case is that there should exist an initial population $x_0$ such
that $\lim_{k \rarr \infty}x_k$ exists and is nonzero. This may be
proved by means of Theorem \ref{PF1}(c).

\msn We now turn to the case of the Leslie model, where $T$ is a
matrix with nonzero elements on the first subdiagonal and $0$'s
elsewhere, and $F$ is a matrix all of whose nonzero elements are
in its first row, see \cite{L1}, \cite[p.38]{Pol} or
\cite[p.4]{C}:
 $$ T = \begin{pmatrix}
    0      & \cdots & 0       &    0  \\
    t_1    & \cdots & 0       &    0  \\
    \vdots & \ddots & \vdots  & \vdots\\
    0      & \cdots &  t_{n-1} & 0
  \end{pmatrix} \qquad \mbox{and} \qquad
     F = \begin{pmatrix}
    f_1   & \cdots  & f_{n-1} &     f_n\\
     0     & \cdots &    0    &      0 \\
    \vdots & \ddots &  \vdots &  \vdots\\
    0      & \cdots &     0   &     0
    \end{pmatrix}.
 $$
 Whether $P$ is irreducible or not, the situation is simpler
as there is only one class of newborns, see \cite{Par} for a
discussion of reducible Leslie models. In this case, the \nrv\
$R_0$ is the $(1,1)$ entry of $Q$ and equals the sum of the
elements of the first row of $F$, as observed in \cite[p.9,l.7]{C}
where it is remarked that the \nrv\ is ``the expected number of
offspring per newborn over the course of its lifetime''. Thus in
this case the \nrv\ equals the net reproductive rate as defined in
the early papers \cite{Ber}, \cite{Lew} and \cite[p.234]{L2}. In
turn, this was an adaptation of a concept previously used in
continuous population models, see \cite[p.115]{Lot}. By Theorem
\ref{compg} we have $R_0 = 1$ if and only if $r = 1$, and if $P$
is a primitive projection matrix for a Leslie model then this has
the interpretation that every population is eventually stationary
if and only if the expected number of offspring over the course of
a newborn's lifetime is $1$, viz. the $(1,1)$ entry of $Q$ is $1$.

\msn Let $P= T+F$ again be a Leslie matrix. Then $q(s)$ is the
leading entry of $F(I-T/s)^{-1}/s$ which is
\begin{equation} \label{poly}
  q(s) = f_1 s^{-1} + f_2 t_1 s^{-2} +
 \cdots + f_n (t_{n-1} \cdots t_1) s^{-n},
\end{equation}
 a polynomial in $s^{-1}$ with coefficients involving all nonzero entries in
 $F$ and $T$. It is classical that the nonzero eigenvalues
$\lambda$ of $P$ satisfy the equation $q(\lambda) = 1$, e.g.
\cite[p.42]{Pol}, and this is consistent with our more general
results since by Theorem \ref{genthm} the growth rate $r$ of $P$
satisfies $q(r) = 1$. We note that the operation of scaling the
fertility matrix by a constant factor was previously considered by
Leslie in \cite[Sec. 5(b)]{L2}. His results in this area have a
somewhat complicated appearance since he wished to determine $s$
as a function of $q(s)$, which requires the solution of a
polynomial equation. This gives rise to the following observation:
To determine the growth rate for a Leslie matrix one needs to find
a positive root of a polynomial equation, but to scale the
fertility matrix to achieve an assigned growth rate $s$ it
suffices to divide the fertility matrix by $q(s)$ in (\ref{poly}).

\msn  Care needs to be taken to interpret the various terms used
in this article when the matrix $P$ is reducible. It is
instructive to consider a simple example. Let
 $$P =  \begin{pmatrix}
          1 & 0 \\
          1 & 1
        \end{pmatrix}. $$
For the projection matrix $P$, an initial population is stable if
and only if its first element is $0$, otherwise it will tend to
infinity. The growth rate $r = \rho(P) = 1$, and the net
reproductive rate $R_0 = 1$ whatever may be the transition matrix
$T$ and fertility matrix $F$ chosen subject to conditions
(\ref{equ1}), (\ref{equ2}) and (\ref{mort}).

\msn We conclude by giving a numerical example which is based on a
plant lifecycle involving vegetative as well as seed reproduction,
\cite[Example 1.c]{Cas1}. Let
\begin{equation*}
 T = (1/2) \begin{pmatrix}
     0 &   0&    0 &   0 &   0 \\
     1 &   0&    0 &   0 &   0 \\
     1 &   0&    0 &   0 &   0 \\
     0 &   1&    1 &   0 &   0  \\
     0  &  0&    0 &   1 &   0
\end{pmatrix} \qquad \mbox{and} \qquad
 F = (1/2) \begin{pmatrix}
     0 &   0&    0 &   0 &   1 \\
     0 &   0&    0 &   0 &   0 \\
     0 &   0&    0 &   1 &   0 \\
     0 &   0&    0 &   0 &   0  \\
     0  &  0&    0 &   0 &   0
\end{pmatrix}.
 \end{equation*}
Let $P = F+T$. This matrix is irreducible, but imprimitive. Then
the growth rate $r = \sqrt 2/2$ and the stable populations for $P$
are $ u = (\sqrt 2, 1, 3, 2\sqrt 2,2)^t$ and its positive
multiples. Further the next generation matrix is
  $$ Q = (1/8)
\begin{pmatrix}
     1 &   1&    1 &   2  &  4 \\
     0 &   0&    0 &   0 &   0 \\
     2 &   2&    2 &   4 &   0 \\
     0 &   0&    0 &   0 &   0  \\
     0 &   0&    0 &   0 &   0
\end{pmatrix}, $$
indicating that the population of newborns that reproduces itself
with the same distribution of newborns in the next generation is
$w = (1, 0, 2, 0, 0)^t$, the Perron vector of the next generation
matrix $Q$. Note that the submatrix in rows and columns $1$ and
$3$ of $Q$ is irreducible and that $R_0 = \rho(Q) = 3/8$. We also
note that $(P/\rho(P))^k w,\ k = 0,1,2, \ldots$ does not tend to a
limit, which indicates that generational stability is compatible
with  permanent oscillation of the normalized population
distribution over time.

\msn If $P$ is replaced by $P_1 = T + 8F/3$, then correspondingly
$Q_1 = (8F/3)(I-T)^{-1} = 8Q/3$ and we have both $\rho(P_1) =
\rho(Q_1) = 1$. Since $\rho(T) = 0$, given any $s > 0$, we can use
(\ref{assign}) to compute $q(s)$ such that $P(s) = T + F/q(s)$ has
growth rate $s$. We obtain $q(s) = (1+2s^2)/8s^4$, and hence the
corresponding net reproductive rate is $R_0(s) = 3s^4/(1+2s^2)$.
The corresponding stable population is $(4s^3, 2s^2, 2s^2 + 8s^4,
2s+4s^3, 1+2s^2)^t$.

\msn {\bf Acknowledgment:} We thank K.P. Hadeler, M.N. Neumann, P.
Ney and several anonymous referees for comments which have helped
to improve this paper.


\end{document}